\documentclass{amsart}
\usepackage{amsmath,amssymb}
\usepackage{yhmath}

\newcommand{\bdis}{\begin{displaymath}}
\newcommand{\edis}{\end{displaymath}}
\newcommand{\be}{\begin{equation}}
\newcommand{\ee}{\end{equation}}
\newcommand{\mbb}{\mathbb}
\newcommand{\mcal}{\mathcal}

\newcommand{\vp}{\varphi}

\newcommand{\zf}{\zeta\left(\frac{1}{2}+it\right)}
 
\newcommand{\zfn}{\zeta\left(\frac{1}{2}+it_\nu\right)}   
\newcommand{\zfnn}{\zeta\left(\frac{1}{2}+it_{\nu+1}\right)}

\DeclareMathOperator{\im}{Im}

\theoremstyle{definition}

\theoremstyle{remark}
\newtheorem{remark}[]{Remark}

\newtheorem*{mydef11}{{\bf Theorem 1}}

\newtheorem*{mydef12}{{\bf Theorem 2}}

\newtheorem*{mydef51}{{\bf Lemma 1}}

\newtheorem*{mydef52}{{\bf Lemma 2}}

\newtheorem*{mydef53}{{\bf Lemma 3}}

\newtheorem*{mydef54}{{\bf Lemma 4}}

\newtheorem*{mydef55}{{\bf Lemma 5}}

\newtheorem*{mydef56}{{\bf Lemma 6}}

\newtheorem*{mydef57}{{\bf Lemma 7}}

\newtheorem*{mydef81}{{\bf Property 1}}

\numberwithin{equation}{section}

\begin{document}

\title[Jacob's ladders, almost linear increments \dots ]{Jacob's ladders, almost linear increments of the Hardy-Littlewood integral (1918) and their relation to the Titchmarsh's sums (1943) and the Fermat-Wiles theorem}

\author{Jan Moser}

\address{Department of Mathematical Analysis and Numerical Mathematics, Comenius University, Mlynska Dolina M105, 842 48 Bratislava, SLOVAKIA}

\email{jan.mozer@fmph.uniba.sk}

\keywords{Riemann zeta-function}

\begin{abstract}
In this paper we give some new consequences that follow from our formula for increments of the Hardy-Littlewood integral. The main of these ones are $\mcal{T}_1$ and $\mcal{T}_2$ equivalents of the Fermat-Wiles theorem.   
\end{abstract}
\maketitle

\section{Introduction} 

\subsection{}    

Let us remind that in our paper \cite{8} we have proved, for example, the following $\zeta$-condition 
\be \label{1.1} 
\lim_{\tau\to \infty}\frac{1}{\tau}\int_{\frac{x^n+y^n}{z^n}\frac{\tau}{1-c}}^{[\frac{x^n+y^n}{z^n}\frac{\tau}{1-c}]^1}\left|\zf\right|^2{\rm d}t\not=1, 
\ee  
where 
\be \label{1.2} 
\left[\frac{x^n+y^n}{z^n}\frac{\tau}{1-c}\right]^1=\vp_1^{-1}\left(\frac{x^n+y^n}{z^n}\frac{\tau}{1-c}\right), 
\ee 
on the class of all Fermat's rationals 
\be \label{1.3} 
\frac{x^n+y^n}{z^n},\ x,y,z\in\mbb{N},\ n\geq 3, 
\ee 
represents the first $\zeta$-equivalent of the Fermat-Wiles theorem, and $\vp_1^{-1}(T)$ denotes the first reverse iteration of the Jacob's ladder $\vp_1(T)$. 

\subsection{} 

Next, let us remind the Dirichlet's function 
\be \label{1.4} 
D(x)=\sum_{n\leq x}d(n), 
\ee 
where $d(n)$ denotes the number of divisors of $n$, and 
\be \label{1.5} 
D(x)=D(N),\ x\in [N,N+1),\ \forall N\in\mbb{N}. 
\ee 
In the paper \cite{9} we have proved, in this direction, the following: the $D$-condition 
\be \label{1.6} 
\lim_{\tau\to\infty}\frac{1}{\tau}\left\{
D\left(\left[\frac{x^n+y^n}{z^n}\frac{\tau}{1-c}\right]^1\right)-D\left(\frac{x^n+y^n}{z^n}\frac{\tau}{1-c}\right)
\right\}\not=1
\ee 
on the class of all Fermat's rationals represents the $D$-equivalent of the Fermat-Wiles theorem. 

\subsection{} 

Further, let us remind: 
\begin{itemize}
	\item[(a)] The Riemann-Siegel formula 
	\be \label{1.7} 
	Z(t)=2\sum_{n\leq\bar{t}}\frac{1}{\sqrt{n}}\cos\{\vartheta(t)-t\ln n\}+\mcal{O}(t^{-1/4}),\ \bar{t}=\sqrt{\frac{t}{2\pi}}, 
	\ee 
	(comp. \cite{11}, p. 79), where 
	\be \label{1.8} 
	\begin{split}
	& \vartheta(t)=-\frac 12t\ln\pi+\im\ln\Gamma\left(\frac{1}{4}+i\frac t2\right)= \\ 
	& \frac{t}{2}\ln\frac{t}{2\pi}-\frac{t}{2}-\frac{\pi}{8}+\mcal{O}\left(\frac 1t\right), 
	\end{split}
	\ee 
	and 
	\be \label{1.9} 
	Z(t)=e^{i\vartheta(t)}\zf. 
	\ee 
	\item[(b)] The Gram's sequence 
	\be \label{1.10} 
	\{ t_\nu\}_{\nu=1}^\infty :\ \vartheta(t_\nu)=\pi\nu. 
	\ee  
	\item[(c)] And, of course, 
	\be \label{1.11} 
	Z(t_\nu)=2(-1)^\nu\sum_{n\leq\bar{t}_\nu}\frac{\cos(t_\nu\ln n)}{\sqrt{n}}+\mcal{O}(t_\nu^{-1/4}),\ \bar{t}_\nu=\sqrt{\frac{\bar{t}_\nu}{2\pi}}. 
	\ee 
\end{itemize} 

\subsection{} 

Finally, let us remind corrected classical Titchmarsh's formulae: 
\be \label{1.12} 
\sum_{\nu=1}^N(-1)^\nu Z(t_\nu)=2N+\mcal{O}(N^{3/4}\ln^{1/4}N), 
\ee 
\be \label{1.13} 
\sum_{\nu=1}^N Z(t_\nu)Z(t_{\nu+1})=-2(1+c)N+\mcal{O})(N^{11/12}\ln N), 
\ee 
(see \cite{10}, p. 101 and p. 105, respectively), where we donote by $c$ the Euler's constant. 

\begin{remark}
Namely, in the formula\footnote{See \cite{10}, (2).} 
\bdis 
f(t)=e^{i\vartheta(t)}\zf=\sum_{n=1}^k\frac{\cos\{\vartheta-t\ln n\}}{\sqrt{n}}+\mcal{O}(t^{-1/4})
\edis  
was omitted the factor $2$ in the last sum, i. e. 
\bdis 
f(t)=\frac 12Z(t). 
\edis 
Moreover, we may put 
\bdis 
\sum_{\nu=M+1}^N(-1)^\nu Z(t_\nu) \ \to \ \sum_{\nu=1}^N (-1)^\nu Z(t_\nu), \dots 
\edis  
into (\ref{1.12}) and (\ref{1.13}), since 
\bdis 
\sum_{\nu=1}^M (-1)^\nu Z(t_\nu)=\mcal{O}(1),\ N\to\infty,\ \dots 
\edis 
\end{remark} 

In this paper we use the following forms 
\be \label{1.14} 
\begin{split}
& \sum_{\nu=1}^N\zfn=2N+\mcal{O}(N^{3/4}\ln^{1/4}N), 
\end{split}
\ee 
\be \label{1.15} 
\sum_{\nu=1}^N\zfn\zfnn=2(1+c)N+\mcal{O}(N^{11/12}\ln N)
\ee  
of the Titchmarsh's formulae (\ref{1.12}), (\ref{1.13}), since\footnote{See (\ref{1.9}), (\ref{1.10}).}
\be \label{1.16} 
\begin{split} 
& (-1)^\nu Z(t_\nu)=\zfn, \\ 
& Z(t_\nu)Z(t_{\nu+1})=-\zfn\zeta\left(\frac 12+it_{\nu+1}\right). 
\end{split} 
\ee  

\subsection{} 
In this paper we obtain, for example, the following result: \\ 

Let 
\be \label{1.17} 
\mcal{T}_1(X)=\sum_{t_\nu\leq X}\zfn. 
\ee  
Then the $\mcal{T}_1$-condition 
\be \label{1.18} 
\lim_{\tau\to\infty}\frac{1}{\tau}\left\{
\mcal{T}_1\left(\left[\frac{x^n+y^n}{z^n}\frac{\tau}{1-c}\right]^1\right)-\mcal{T}_1\left(\frac{x^n+y^n}{z^n}\frac{\tau}{1-c}\right)
\right\}\not=\frac{1}{\pi}
\ee 
on the class of all Fermat's rationals represents the $\mcal{T}_1$-equivalent of the Fermat-Wiles theorem.  

Next, let 
\be \label{1.19} 
\mcal{T}_2(X)=\sum_{t_\nu\leq x}\zfn\zfnn. 
\ee  
Then the $\mcal{T}_2$-condition 
\be \label{1.20} 
\lim_{\tau\to\infty}\frac{1}{\tau}\left\{
\mcal{T}_2\left(\left[\frac{x^n+y^n}{z^n}\frac{\tau}{1-c}\right]^1\right)-\mcal{T}_2\left(\frac{x^n+y^n}{z^n}\frac{\tau}{1-c}\right)
\right\}\not=\frac{1+c}{\pi}
\ee 
on the class of all Fermat's rationals represents the $\mcal{T}_2$-equivalent of the Fermat-Wiles theorem. 

\subsection{} 

In this paper we use the following notions from our papers \cite{2} -- \cite{5}: 
\begin{itemize}
	\item[{\tt (a)}] Jacob's ladder $\vp_1(t)$, 
	\item[{\tt (b)}] the function 
	\bdis 
	\begin{split}
		& \tilde{Z}^2(t)=\frac{{\rm d}\vp_1(t)}{{\rm d}t}=\frac{1}{\omega(t)}\left|\zf\right|^2,\\ 
		& \omega(t)=\left\{1+\mcal{O}\left(\frac{\ln\ln t}{\ln t}\right)\right\}\ln t,\ t\to\infty, 
	\end{split}
	\edis 
	\item[{\tt (c)}] direct iterations of Jacob's ladders 
	\bdis 
	\begin{split}
		& \vp_1^0(t)=t,\ \vp_1^1(t)=\vp_1(t),\ \vp_1^2(t)=\vp_1(\vp_1(t)),\dots , \\ 
		& \vp_1^k(t)=\vp_1(\vp_1^{k-1}(t))
	\end{split}
	\edis 
	for every fixed natural number $k$, 
	\item[{\tt (d)}] reverse iterations of Jacob's ladders 
	\bdis 
	\begin{split}
		& \vp_1^{-1}(T)=\overset{1}{T},\ \vp_1^{-2}(T)=\vp_1^{-1}(\overset{1}{T})=\overset{2}{T},\dots, \\ 
		& \vp_1^{-r}(T)=\vp_1^{-1}(\overset{r-1}{T})=\overset{r}{T},\ r=1,\dots,k, 
	\end{split} 
	\edis 
	where, for example, 
	\be \label{1.21} 
	\vp_1(\overset{r}{T})=\overset{r-1}{T}
	\ee  
	for every fixed $k\in\mbb{N}$, and  
	\be \label{1.22}
	\begin{split} 
		& \overset{r}{T}-\overset{r-1}{T}\sim(1-c)\pi(\overset{r}{T});\ \pi(\overset{r}{T})\sim\frac{\overset{r}{T}}{\ln \overset{r}{T}},\ r=1,\dots,k,\ T\to\infty, \\ 
		& \overset{0}{T}=T<\overset{1}{T}(T)<\overset{2}{T}(T)<\dots<\overset{k}{T}(T), \\ 
		& T\sim \overset{1}{T}\sim \overset{2}{T}\sim \dots\sim \overset{k}{T},\ T\to\infty. 
	\end{split}
	\ee 
\end{itemize} 

\begin{remark}
	The asymptotic behaviour of the points 
	\bdis 
	\{T,\overset{1}{T},\dots,\overset{k}{T}\}
	\edis  
	is as follows: at $T\to\infty$ these points recede unboundedly each from other and all together are receding to infinity. Hence, the set of these points behaves at $T\to\infty$ as one-dimensional Friedmann-Hubble expanding Universe. 
\end{remark}

\section{Jacob's ladders and \emph{proliferation} of every $L_2$-orthogonal system} 

Let us remind that we may view the results (\ref{1.1}), (\ref{1.6}), (\ref{1.18}) and (\ref{1.20}) as points of contact between functions 
\bdis 
\left|\zf\right|^2, D(x), \mcal{T}_1(X), \mcal{T}_2(X)
\edis 
and the Fermat-Wiles theorem. 

\begin{remark}
The basis for above mentioned phenomena is constituted by Jacob's ladders, see \cite{2}, and almost linear increments of the classical Hardy-Littlewood integral, see \cite{7}. 
\end{remark} 

For completeness we remind also our result, see \cite{6}, that gives the point of contact between the function 
\bdis 
\left|\zf\right|^2
\edis 
and the theory of $L_2$-orthogonal systems. 

\subsection{} 

We have introduced the generating vector operator $\hat{G}$ acting on the class of all $L_2$-orthogonal systems 
\be \label{2.1} 
\{f_n(t)\}_{n=0}^\infty,\ t\in [a,a+2l],\ a\in\mbb{R},\ l>0
\ee 
as 
\be \label{2.2} 
\begin{split}
	& \{ f_n(t)\}_{n=0}^\infty\xrightarrow{\hat{G}}\{ f^{p_1}_n(t)\}_{n=0}^\infty\xrightarrow{\hat{G}}\{ f^{p_1,p_2}_n(t)\}_{n=0}^\infty\xrightarrow{\hat{G}} \dots \\ 
	& \xrightarrow{\hat{G}}\{ f^{p_1,p_2,\dots,p_s}_n(t)\}_{n=0}^\infty,\ p_1,\dots,p_s=1,\dots,k 
\end{split}
\ee  
for every fixed $k,s\in\mbb{N}$ with explicit formulae\footnote{See \cite{6}, (2.19).} for 
\bdis 
f^{p_1,p_2,\dots,p_s}_n(t). 
\edis 

\subsection{} 

In the case of Legendre's orthogonal system 
\be \label{2.3} 
\{ P_n(t)\}_{n=0}^\infty,\ t\in [-1,1]
\ee 
the operator $\hat{G}$ produces, for example, the third generation as follows 
\be \label{2.4} 
\begin{split}
	& P_n^{p_1,p_2,p_3}(t)=P_n(u_{p_1}(u_{p_2}(u_{p_3}(t))))\times\prod_{r=0}^{p_1-1}\left|\tilde{Z}(v_{p_1}^r(u_{p_2}(u_{p_3}(t))))\right|\times \\ 
	& \prod_{r=0}^{p_2-1}\left|\tilde{Z}(v_{p_2}^r(u_{p_3}(t)))\right|\times 
	\prod_{r=0}^{p_3-1}\left|\tilde{Z}(v_{p_3}^r(t))\right|, \\ 
	& p_1,p_2,p_3=1, \dots , k,\ t\in [-1,1],\ a=-1,\ l=1, 
\end{split}
\ee 
where 
\be \label{2.5} 
u_{p_i}(t)=\vp_1^{p_i}\left(\frac{\overset{p_i}{\wideparen{T+2}}-\overset{p_i}{T}}{2}(t+1)-\overset{p_i}{T}\right)-T-1,\ i=1,2,3
\ee 
are automorphisms on $[-1,1]$ and 
\be \label{2.6} 
\begin{split}
	& v_{p_i}^r(t)=\vp_1^r\left(\frac{\overset{p_i}{\wideparen{T+2}}-\overset{p_i}{T}}{2}(t+1)-\overset{p_i}{T}\right),\ r=0,1,\dots,p_i-1, \\ 
	& t\in[-1,1]\Rightarrow u_{p_i}(t)\in [-1,1]\wedge v_{p_i}^r(t)\in[\overset{p_i-r}{T},\overset{p_i-r}{\wideparen{T+2}}]. 
\end{split}
\ee 

\begin{mydef81}
	\begin{itemize}
		\item[(a)] Every member of every new $L_2$-orthogonal system 
		\be \label{2.7} 
		\{P_n^{p_1,p_2,p_3}(t)\}_{n=0}^\infty,\ t\in[-1,1],\ p_1,p_2,p_3=1,\dots,k
		\ee 
		contains the function 
		\bdis 
		\left|\zf\right|_{t=\tau}
		\edis 
		for corresponding $\tau$ since\footnote{See \cite{3}, (9.1), (9.2).} 
		\be \label{2.8} 
		|\tilde{Z}(t)|=\sqrt{\frac{{\rm d}\vp_1(t)}{{\rm d}t}}=\frac{\{1+o(1)\}}{\sqrt{\ln t}}\left|\zf\right|,\ t\to\infty; 	
		\ee 
		\item[(b)] Property (a) holds true due to Theorem of the paper \cite{6} for every generation 
		\bdis 
		\{f_n^{p_1,\dots,p_s}(t)\}_{n=0}^\infty,\ t\in [a,a+2l],\ s\in\mbb{N}. 
		\edis 
	\end{itemize}
\end{mydef81} 

\begin{remark}
	Our type of \emph{proliferation} of every $L_2$-orthogonal system is in context with the Chumash, Bereishis, 26:12, \emph{Isaac sowed in the land, and in that year reaped a hundredfold, thus had HASHEM blessed him}. 
\end{remark} 

\subsection{} 

Next, according to (\ref{2.4}), to the $L_2$-orthogonal system 
\be \label{2.9} 
\{P_n^{p_1,p_2,p_3}(t)\}_{n=0}^\infty 
\ee  
corresponds, for example, the following $L_2$-orthonormal system 
\be \label{2.10} 
\{\bar{P}_n^{p_1,p_2,p_3}(t)\}_{n=0}^\infty,\ t\in [-1,1],\ a=-1,\ l=1, 
\ee  
where\footnote{Comp. \cite{6}, (2.20).}  
\be \label{2.11} 
\bar{P}_n^{p_1,p_2,p_3}(t)=\left(\prod_{i=0}^3\sqrt{\frac{2}{\overset{i}{\wideparen{T+2}}-\overset{i}{T}}}\right)P_n^{p_1,p_2,p_3}(t). 
\ee  

\subsection{} 

Now, we have the following property as the consequence of the Menshow-Rademacher theorem: If the sequence 
\be \label{2.12} 
\{a_n^{p_1,p_2,p_3}\}_{n=0}^\infty,\ a_n^{p_1,p_2,p_3}\in\mbb{R} 
\ee 
fulfils the condition 
\be \label{2.13} 
\sum_{n=0}^\infty \{a_n^{p_1,p_2,p_3}\ln(n+1)\}^2<+\infty, 
\ee  
then the orthogonal series 
\be \label{2.14} 
\sum_{n=0}^\infty a_n^{p_1,p_2,p_3}\bar{P}_n^{p_1,p_2,p_3}(t) 
\ee  
converges almost everywhere on $[-1,1]$, i. e. there is a function 
\be \label{2.15} 
\begin{split}
& F^{p_1,p_2,p_3}(t)= F^{p_1,p_2,p_3}(t;\{\bar{P}_n^{p_1,p_2,p_3}(t)\}_{n=0}^\infty,\{a_n^{p_1,p_2,p_3}\}_{n=0}^\infty)
\end{split}
\ee 
such that the equality 
\be \label{2.16} 
F^{p_1,p_2,p_3}(t)=\sum_{n=0}^\infty a_n^{p_1,p_2,p_3}\bar{P}_n^{p_1,p_2,p_3}(t) 
\ee 
holds true almost everywhere on $[-1,1]$. 

\begin{remark}
Of course, for every generation 
\bdis 
\{f_n^{p_1,\dots,p_s}(t)\}_{n=0}^\infty,\ t\in [a,a+2l],\ s\in\mbb{N},\ l\in\mbb{R}^+, 
\edis 
there are analogues of the formulae (\ref{2.9})--(\ref{2.16}). 
\end{remark} 

\section{Variant of the first Titchmarsh's formula} 

\subsection{} 

Since, see (\ref{1.10}) 
\be \label{3.1} 
\vartheta(t_N)=\pi N, 
\ee  
we have by (\ref{1.8}) 
\be \label{3.2} 
\pi N=\frac{t_N}{2}\ln\frac{t_N}{2\pi}-\frac{t_N}{2}-\frac{\pi}{8}+\mcal{O}\left(\frac{1}{t_N}\right). 
\ee 
Consequently we obtain 
\be \label{3.3} 
\begin{split}
& 2N=\frac{1}{\pi}t_N\ln\frac{t_N}{2\pi}-\frac{t_N}{\pi}-\frac{1}{4}+\mcal{O}\left(\frac{1}{t_N}\right)= \\ 
& \frac{1}{\pi}t_N\ln t_N-\frac{1}{\pi}(1+\ln 2\pi)t_N+\mcal{O}(1). 
\end{split}
\ee 

\subsection{} 

Now, we put 
\be \label{3.4} 
\mcal{T}_1(t_N)=\sum_{\nu\leq t_N}\zfn. 
\ee  
Then we obtain, by making use of (\ref{3.3}) in the formula (\ref{1.14}), the following formula 
\be \label{3.5} 
\mcal{T}_1(t_N)=\frac{1}{\pi}t_N\ln t_N-\frac{1}{\pi}(1+\ln 2\pi)t_N+\mcal{O}(t_N^{3/4}\ln t_N), 
\ee  
since, of course, 
\be \label{3.6} 
\mcal{O}(N^{3/4}\ln^{1/4}N)=\mcal{O}(t_N^{3/4}\ln t_N). 
\ee  

\subsection{} 

Next, let us remind that\footnote{See \cite{10}, p. 102.}
\be \label{3.7} 
t_{N+1}-t_N\sim \frac{2\pi}{\ln t_N},\ N\to\infty 
\ee  
and\footnote{Comp. \cite{11}, p. 99.} 
\be \label{3.8} 
\zf=\mcal{O}(t^{1/6}),\ t\to\infty. 
\ee 
Since the error term produced by the substitution 
\be \label{3.9} 
t_N\to x,\ x\in [t_N,t_{N+1}),\ \forall N\in\mbb{N}
\ee 
is comfortably absorbed into $\mcal{O}(t_N^{3/4}\ln t_N)$ by (\ref{3.7}) and (\ref{3.8}), then we obtain the following statement. 

\begin{mydef51}
If 
\be \label{3.10} 
\mcal{T}_1(X)=\sum_{t_\nu\leq x}\zfn, 
\ee  
where 
\be \label{3.11} 
\mcal{T}_1(X)=\mcal{T}_1(t_N),\ \forall X\in [t_N,t_{N+1}),\ \forall N\in\mbb{N}, 
\ee  
then we obtain the following variant of the first Titchmarsh's formula 
\be \label{3.12} 
\mcal{T}_1(X)=\frac{1}{\pi}X\ln X-\frac{1}{\pi}(1+\ln 2\pi)X+\mcal{O}(X^{3/4}\ln X),\ X\to\infty. 
\ee 
\end{mydef51} 

\section{Connection between increments of the Hardy-Littlewood integral and increments of the first Titchmarsh's function $\mcal{T}_1(x)$} 

\subsection{} 

Now we will continue by the similar way as in the case of elementary Dirichlet's formula for the sum of number of divisors in our paper \cite{9}. We will use: 
\begin{itemize}
	\item[(a)] Our almost exact formula, \cite{9}, (3.7), in the case 
	\be \label{4.1} 
	T\to \overset{r}{X},\ r=1,\dots,k, 
	\ee  
	where 
	\be \label{4.2} 
	\overset{r}{X}=\vp_1^t(X) 
	\ee  
	for every fixed natural number $k$, i. e. the formula ($\overset{r}{X}\sim X$) 
	\be \label{4.3} 
	\begin{split}
	& \frac{1}{\pi}\int_0^{\overset{r}{X}}\left|\zf\right|^2{\rm d}t=\frac{1}{\pi}\overset{r-1}{X}\ln\overset{r-1}{X}+ \\ 
	& \frac{1}{\pi}(c-\ln 2\pi)\overset{r-1}{X}+c_0+\mcal{O}\left(\frac{\ln X}{X}\right). 
	\end{split}
	\ee 
	\item[(b)]  Our formula (\ref{3.12}) in the case 
	\be \label{4.4} 
	T\to\overset{r-1}{X},\ \overset{0}{X}=X, 
	\ee 
	that implies 
	\be \label{4.5} 
	\mcal{T}_1(\overset{r-1}{X})=\frac{1}{\pi}\overset{r-1}{X}\ln\overset{r-1}{X}-\frac{1}{\pi}(1+\ln 2\pi)\overset{r-1}{X}+\mcal{O}(X^{3/4}\ln X). 
	\ee 
\end{itemize} 

Next, by subtracting (\ref{4.5}) from (\ref{4.3}) we obtain 
\be \label{4.6} 
\begin{split}
& \frac{1}{\pi}\int_0^{\overset{r}{X}}\left|\zf\right|^2{\rm d}t-\mcal{T}_1(\overset{r-1}{X})= \\ 
& \frac{1+c}{\pi}\overset{r-1}{X}+\mcal{O}(X^{3/4}\ln X), 
\end{split}
\ee 
and the translation $r\to r+1$ gives us 
\be \label{4.7} 
\frac{1}{\pi}\int_0^{\overset{r+1}{X}}\left|\zf\right|^2{\rm d}t-\mcal{T}_1(\overset{r}{X})=\frac{1+c}{\pi}\overset{r}{X}+\mcal{O}(X^{3/4}\ln X). 
\ee 
And consequently, subtraction of (\ref{4.6}) from (\ref{4.7}) gives the following formula 
\be \label{4.8} 
\begin{split}
& \frac{1}{\pi}\int_{\overset{r}{X}}^{\overset{r+1}{X}}\left|\zf\right|^2{\rm d}t-\mcal{T}_1(\overset{r}{X})+\mcal{T}_1(\overset{r-1}{X})= \\ 
& \frac{1+c}{\pi}(\overset{r}{X}-\overset{r-1}{X})+\mcal{O}(X^{3/4}\ln X),\ r=1,\dots,k. 
\end{split}
\ee 
Next, if we use the formula\footnote{This one follows from our almost linear formula, see \cite{7}, with $\delta$ positive and small.} 
\be \label{4.9} 
\begin{split}
& \frac{1}{\pi}\int_{\overset{r}{X}}^{\overset{r+1}{X}}\left|\zf\right|^2{\rm d}t-\frac{1}{\pi}\int_{\overset{r-1}{X}}^{\overset{r}{X}}\left|\zf\right|^2{\rm d}t= \\ 
& (1-c)(\overset{r}{X}-\overset{r-1}{X})+\mcal{O}(X^{1/3+\delta})
\end{split}
\ee 
in (\ref{4.8}), then we obtain 
\be \label{4.10} 
\begin{split}
& \frac{1}{\pi}\int_{\overset{r-1}{X}}^{\overset{r}{X}}\left|\zf\right|^2{\rm d}t=\mcal{T}_1(\overset{r}{X})-\mcal{T}_1(\overset{r-1}{X})+ \\ 
& \frac{2c}{\pi}(\overset{r}{X}-\overset{r-1}{X})+\mcal{O}(X^{3/4}\ln X). 
\end{split}
\ee 
Since\footnote{See (\ref{1.22}).} 
\be \label{4.11} 
\overset{r}{X}-\overset{r-1}{X}\sim (1-c)\pi(\overset{r}{X})\sim(1-c)\frac{\overset{r}{X}}{\ln \overset{r}{X}}\sim (1-c)\frac{X}{\ln X},\ X\to\infty, 
\ee 
we obtain the following lemma. 

\begin{mydef52}
\be \label{4.12} 
\begin{split}
& \mcal{T}_1(\overset{r}{X})-\mcal{T}_1(\overset{r-1}{X})=\frac{1}{\pi}\int_{\overset{r-1}{X}}^{\overset{r}{X}}\left|\zf\right|^2{\rm d}t+\mcal{O}\left(\frac{X}{\ln X}\right), \\ 
& r=1,\dots,k,\ X\to\infty, 
\end{split}
\ee 
comp. with \cite{9}, (4.9). 
\end{mydef52}

\section{$\mcal{T}_1$ and $\mcal{T}_2$ equivalents of the Fermat-Wiles theorem} 

\subsection{} 

In what follows we shall use, for example, the formula (\ref{4.12}), $r=1$, i. e. this one: 
\be \label{5.1} 
\begin{split}
& \mcal{T}_1(\overset{1}{X})-\mcal{T}_1(X)=\frac{1}{\pi}\int_{X}^{\overset{1}{X}}\left|\zf\right|^2{\rm d}t+\mcal{O}\left(\frac{X}{\ln X}\right), \\ 
& X=\overset{0}{X},\ X>X_0>0, 
\end{split}
\ee 
where $x_0$ is sufficiently big and 
\be \label{5.2} 
\overset{1}{X}=[X]^1=\vp_1^{-1}(X). 
\ee 
Now, if we put 
\be \label{5.3} 
X=\frac{x}{1-c}\tau,\ \tau\in \left(\frac{1-c}{x}x_0,+\infty\right),\ x>0
\ee 
into (\ref{5.1}), we obtain the following statement. 

\begin{mydef53}
\be \label{5.4} 
\begin{split}
& \mcal{T}_1\left(\left[\frac{x}{1-c}\tau\right]^1\right)-\mcal{T}_1\left(\frac{x}{1-c}\tau\right)= \\ 
& \frac{1}{\pi}\int_{\frac{x}{1-c}\tau}^{\left[\frac{x}{1-c}\tau\right]^1}\left|\zf\right|^2{\rm d}t+\mcal{O}\left(\frac{\tau}{\ln\tau}\right), \\ 
& \tau\in (\tau_1(x),+\infty),\ \tau_1(x)=\max\left\{\frac{(1-c)^2}{x^2},(x_0)^2\right\}
\end{split}
\ee  
for every fixed $x>0$, where, of course, 
\be \label{5.5} 
\frac{1-c}{x}x_0\leq \tau_1(x),\ x>0, 
\ee 
and the constant in the $\mcal{O}$-term depends on $x$. 
\end{mydef53}  

Next, since it is true that\footnote{See \cite{8}, (4.6).} 
\be \label{5.6} 
\lim_{\tau\to\infty}\frac{1}{\tau}\int_{\frac{x}{1-c}\tau}^{\left[\frac{x}{1-c}\tau\right]^1}\left|\zf\right|^2{\rm d}t=x, 
\ee  
then next lemma follows from (\ref{5.4}). 

\begin{mydef54} 
\be \label{5.7} 
\lim_{\tau\to\infty}\frac{1}{\tau}\left\{\mcal{T}_1\left(\left[\frac{x}{1-c}\tau\right]^1\right)-\mcal{T}_1\left(\frac{x}{1-c}\tau\right)\right\}=\frac{1}{\pi}x
\ee 
for every fixed $x>0$, where\footnote{See (\ref{5.2}).}  
\be \label{5.8} 
\left[\frac{x}{1-c}\tau\right]^1=\vp_1^{-1}\left(\frac{x}{1-c}\tau\right). 
\ee 
\end{mydef54}   

\subsection{} 

Now, if we use the substitution 
\be \label{5.9} 
x\to \frac{x^n+y^n}{z^n},\ x,y,z\in\mbb{N},\ n\geq 3
\ee 
in (\ref{5.7}), then we obtain the following. 

\begin{mydef55}
\be \label{5.10} 
\begin{split}
& \lim_{\tau\to\infty}\frac{1}{\tau}\left\{\mcal{T}_1\left(\left[\frac{x^n+y^n}{z^n}\frac{\tau}{1-c}\right]^1\right)-\mcal{T}_1\left(\frac{x^n+y^n}{z^n}\frac{\tau}{1-c}\right)\right\}= \\ 
& \frac{1}{\pi}\frac{x^n+y^n}{z^n}
\end{split}
\ee 
for every fixed Fermat's rational 
\bdis 
\frac{x^n+y^n}{z^n}. 
\edis 
\end{mydef55} 

Consequently, we have the following theorem. 

\begin{mydef11}
The $\mcal{T}_1$-condition 
\be \label{5.11} 
\lim_{\tau\to\infty}\frac{1}{\tau}\left\{\mcal{T}_1\left(\left[\frac{x^n+y^n}{z^n}\frac{\tau}{1-c}\right]^1\right)-\mcal{T}_1\left(\frac{x^n+y^n}{z^n}\frac{\tau}{1-c}\right)\right\}\not=\frac{1}{\pi}
\ee 
on the class of all Fermat's rationals represents the $\mcal{T}_1$-equivalent of the Fermat-Wiles theorem. 
\end{mydef11}

\subsection{} 

Next, it is clear that our method (\ref{3.1})--(\ref{3.12}) gives the following result for the second Titchmarsh's formula (\ref{1.15}). 

\begin{mydef56}
\be \label{5.12} 
\begin{split}
& \mcal{T}_2(X)=\sum_{t_\nu\leq X}\zfn\zfnn= \\ 
& \frac{1+c}{\pi}X\ln X-\frac{1+c}{\pi}(1+\ln 2\pi)X+\mcal{O}(X^{11/12}\ln^{23/12}X), 
\end{split}
\ee  
where 
\be \label{5.13} 
\mcal{T}_2(X)=\mcal{T}_2(t_N),\ \forall X\in [t_N,t_{N+1}), \ \forall N\in\mbb{N}.
\ee 
\end{mydef56} 

Hence, we obtain by methods (\ref{4.1})--(\ref{4.12}) and (\ref{5.1})--(\ref{5.10}) next lemma. 
\begin{mydef57}
The following formula holds true for every fixed Fermat's rational 
\be \label{5.14} 
\begin{split}
& \lim_{\tau\to\infty}\frac{1}{\tau}\left\{\mcal{T}_2\left(\left[\frac{x^n+y^n}{z^n}\frac{\tau}{1-c}\right]^1\right)-\mcal{T}_2\left(\frac{x^n+y^n}{z^n}\frac{\tau}{1-c}\right)\right\}= \\ 
& \frac{1+c}{\pi}\frac{x^n+y^n}{z^n}. 
\end{split}
\ee 
\end{mydef57} 

And finally, our second theorem follows. 

\begin{mydef12}
The $\mcal{T}_2$-condition 
\be \label{5.15} 
\lim_{\tau\to\infty}\frac{1}{\tau}\left\{\mcal{T}_2\left(\left[\frac{x^n+y^n}{z^n}\frac{\tau}{1-c}\right]^1\right)-\mcal{T}_2\left(\frac{x^n+y^n}{z^n}\frac{\tau}{1-c}\right)\right\}\not=\frac{1+c}{\pi}
\ee 
on the class of all Fermat's rationals represents the $\mcal{T}_2$-equivalent of the Fermat-Wiles theorem. 
\end{mydef12}

I would like to thank Michal Demetrian for his moral support of my study of Jacob's ladders.

\end{document}